\documentclass[twoside,reqno,A4]{amsart}

\usepackage{latexsym,rotate,eucal,cite,url}
\usepackage{amsmath,amsthm,amssymb,amsxtra} 

\newcommand{\ang}[1]{\langle{#1}\rangle}

\newcommand{\mb}[1]{\mathbb{#1}}

\newcommand{\binm}[2]{\bigg(\genfrac{}{}{0mm}{0}{#1}{#2}\bigg)}
\newcommand{\binq}[2]{\bigg[\genfrac{}{}{0mm}{0}{#1}{#2}\bigg]}
\newcommand{\bing}[2]{\bigg\{\genfrac{}{}{0mm}{0}{#1}{#2}\bigg\}}

\newcommand{\be}{\begin{equation}}
\newcommand{\ee}{\end{equation}}
\newcommand{\ba}{\begin{array}}
\newcommand{\ea}{\end{array}}
\newcommand{\bmn}{\begin{eqnarray}}
\newcommand{\emn}{\end{eqnarray}}
\newcommand{\bnm}{\begin{eqnarray*}}
\newcommand{\enm}{\end{eqnarray*}}
\newcommand{\bln}{\begin{subequations}}
\newcommand{\eln}{\end{subequations}}

\newcommand{\pq}[1]{\begin{equation}#1\end{equation}}
\newcommand{\pnq}[1]{\begin{align*}#1
            \end{align*}}    
\newcommand{\centro}[1]
           {\begin{center}#1\end{center}}

\newcommand{\lam}{\lambda}

\newcommand{\del}{\delta}

\newtheorem{thm}{Theorem}
\newtheorem{lemm}[thm]{Lemma}
\newtheorem{corl}[thm]{Corollary}
\newtheorem{prop}[thm]{Proposition}

\newtheorem{entry}{Entry}

\newcommand{\thank}[2]{\begin{center}\parbox{#1}
{{\sc\bf Acknowledgement}:\,{\small\it#2}}\end{center}}

\newcommand{\referxy}[4]{\bibitem{kn:#1}{#2,}~\emph{#3,}~{#4.}}	
\newcommand{\cito}[1]{\cite{kn:#1}}	
\newcommand{\citu}[2]{\cite[#2]{kn:#1}}

\begin{document}
\title{Convolution Identities of Stirling Numbers}
\author{Nadia N. Li and Wenchang Chu}

\address{\textbf{N.~N.~Li}:
School of Mathematics and Statistics\newline
Zhoukou Normal University, Zhoukou (Henan), China}
\email{lina20190606@outlook.com}
\thanks{*~Corresponding author (N.~N.~Li): lina20190606@outlook.com}
\address{\textbf{W.~Chu}:
Department of Mathematics and Physics\newline
University of Salento, Lecce 73100, Italy}
\email{chu.wenchang@unisalento.it}


\begin{abstract}\vspace{-5mm}
By means of the generating function method, a linear recurrence 
relation is explicitly resolved. The solution is expressed in 
terms of the Stirling numbers of both the first and the second
kind. Two remarkable pairs of combinatorial identities are 
established as applications, that contain some well--known 
convolution formulae on Stirling numbers as special cases.
\end{abstract}

\maketitle\thispagestyle{empty}\vspace*{-5mm} 
\markboth{Nadia N. Li and Wenchang Chu}
{Convolution Identities of Stirling Numbers}

\section{Introduction and Motivation}

Denote by $\mb{N}$ the set of natural numbers with $\mb{N}_0=\mb{N}\cup\{0\}$.
For an indeterminate $x$, define the rising and falling factorials by 
$(x)_0=\ang{x}_0\equiv1$ and  
\pnq{(x)_n&=x(x+1)\cdots(x+n-1)
\quad\text{for}\quad n\in\mb{N},\\
\ang{x}_n&=x(x-1)\cdots(x-n+1)
\quad\text{for}\quad n\in\mb{N}.}

Then the unsigned Stirling numbers of the first kind $\binq{n}{k}$
are determined as the connection coefficients:
\[\ang{x}_n=\sum_{k=0}^n(-1)^{n-k}\binq{n}{k}x^k.\]
Analogously, the Stirling numbers of the second kind $\bing{n}{k}$
are given by
\[x^n=\sum_{k=0}^n\bing{n}{k}\ang{x}_k\]
which admits the explicit expression 
\[\bing{n}{k}=\frac{1}{k!}
\sum_{j=0}^k(-1)^{k-j}\binm{k}{j}j^n.\]
These numbers appear frequently in mathematical literatures 
and have wide applications in combinatorics and number theory
(see \citu{comtet}{Chapter~5}, \citu{graham}{\S6.2} and \cito{knuth},
just for example). 

Recently, Stenlund~\cito{stenl} introduced an interesting  bivariate polynomial 
sequence $\{P_n(x,z)\}_{n\in\mb{N}}$ through the following recurrence relation:
\pq{\label{P-rr}
P_{n+1}(x,z)=x\binm{n+z}{n}
-x\sum_{m=1}^n\binm{n-m+z}{n-m+1}P_m(x,z)
\quad\text{with\quad}P_1(x,z)=x.}
The same author not only found out an explicit double sum expression
\pq{\label{P-double}
P_n(x,z)=\sum_{k=1}^n\sum_{j=1}^k
\frac{(-1)^{j-1}(j-1)!}{(n-1)!}
\binq{n}{k}\bing{k}{j}x^jz^{k-1},}
but also explored applications to combinatorial identities and probability theory.

Inspired by the above work of Stenlund~\cito{stenl},
we shall examine the extended polynomial sequence 
$\{Q_n\}_{n\in\mb{N}_0}$ defined by the recurrence relation
\pq{\label{Q-rr}
Q_{n}=\binm{\lam}{n}(-1)^n
-x\sum_{k=1}^n(-1)^k
\binm{y}{k}Q_{n-k}
\quad\text{with\quad}Q_0=1.}
The first three terms are recorded as follows:
\pnq{Q_1&=xy-\lam,\qquad
Q_2=\frac{x y}{2}(1-y+2 x y)-\frac{\lam}{2}(1-\lam +2 x y),\\
Q_3&=\frac{xy}{6}(y-1)(y-2) -\frac{\lam}{6}(\lam -1)(\lam -2)
-x y (1-y+xy)(\lam -x y)+\frac{\lam  x y}{2}(\lam -y).}

These $Q_n$-polynomials in \eqref{Q-rr} generalize the $P_n$-polynomials 
in \eqref{P-rr} because it is not hard to verify that $Q_n$ becomes
$P_{n+1}/x$ when $\lam=-1-z$ and $y\to -z$.  

The rest of the paper will be organized as follows. In the next section,
we shall derive the generating function and the explicit formulae for
the $P_n$-polynomials. As applications, two pais of convolution sums containing 
both kinds of Stirling numbers will be evaluated in closed forms in Section~3.

\section{Generating Function and Explicit Expressions}
For the sequence $\{Q_n\}_{n\in\mb{N}_0}$ defined by \eqref{Q-rr},
we shall derive, by making use of generating function approach, 
their generating function and explicit formulae. The informed 
reader will notice that this treatment is more transparent than 
the induction approach adopted by Stenlund~\cito{stenl}. 

Multiplying both sides of \eqref{Q-rr} by $\tau^{n}$ and then summing 
over $n$ for $n\geq0$, we can manipulate the generating function
\pnq{Q(\tau):=&\sum_{n=0}^{\infty}
    Q_{n}\tau^{n}
    =1+\sum_{n=1}^{\infty}
    Q_{n}\tau^{n}\\
=&1+\sum_{n=1}^{\infty}
\binm{\lam}{n}(-\tau)^n
-x\sum_{n=1}^{\infty}
\sum_{k=1}^n(-1)^k
\binm{y}{k}Q_{n-k}\tau^n\\
=&(1-\tau)^\lam
-x\sum_{k=1}^{\infty}
(-1)^k\binm{y}{k}\tau^k
\sum_{n=k}^\infty Q_{n-k}\tau^{n-k}.}
This leads us to the functional equation
\[Q(\tau)
=(1-\tau)^\lam
-xQ(\tau)\{(1-\tau)^y-1\}\]
and the rational function expression.
\begin{lemm}[Generating function]\label{Q-gf}
\[Q(\tau)=\frac{(1-\tau)^{\lam}}{1-x+x(1-\tau)^{y}}.\]
\end{lemm}

By writing in the geometric series
\pnq{Q(\tau)&=\frac{1}{1-x}\times
\frac{(1-\tau)^{\lam}}{1-\frac{x}{x-1}(1-\tau)^{y}}\\
&=\sum_{k=0}^\infty
\frac{(-x)^k}{(1-x)^{k+1}}
(1-\tau)^{\lam+ky}}
and then extracting the coefficient of $\tau^{n}$,
we find the infinite series expression.
\begin{prop}[Single sum formula]\label{Q1s}
\[Q_n=\sum_{k=0}^\infty(-1)^{n+k}
\binm{\lam+ky}{n}\frac{x^k}{(1-x)^{k+1}}.\]
\end{prop}

Expanding the rightmost fraction into binomial series,
we have further
\pnq{Q_n&=\sum_{k=0}^\infty
(-1)^{n+k}\binm{\lam+ky}{n}
\sum_{j=0}^\infty
\binm{k+j}{j}x^{k+j}\\
&=\sum_{m=0}^\infty(-1)^{n}x^m\sum_{k=0}^m
(-1)^{k}\binm{m}{k}\binm{\lam+ky}{n},}
where the last line is justified by the replacement $k+j=m$.
Observe that the inner sum with respect to $k$ results substantially
in the differences of order $m$ about a polynomial of degree $n$,
which is equal to zero when $m>n$. Therefore, we have derived 
a finite double sum expression.
\begin{prop}[Double sum formula]\label{Q2s}
\[Q_n=\sum_{m=0}^n\sum_{k=0}^m
(-1)^{n+k}\binm{m}{k}\binm{\lam+ky}{n}x^m.\]
\end{prop}

Expressing the last binomial coefficient in terms 
of unsigned Stirling numbers and then making use 
of the binomial theorem, we have   
\pnq{\binm{\lam+ky}{n}
&=\frac{\ang{\lam+ky}_n}{n!}
=\sum_{i=0}^n\frac{(-1)^{n-i}}{n!}
\binq{n}{i}(\lam+yk)^{i}\\
&=\sum_{i=0}^n\frac{(-1)^{n-i}}{n!}
\binq{n}{i}
\sum_{j=0}^i\binm{i}{j}   
\lam^{i-j}(ky)^j.}

By substitutions, $Q_n$ can be written as a four--fold sum
\pnq{Q_n&=\sum_{m=0}^n\sum_{k=0}^m
\binm{m}{k}x^m
\sum_{i=0}^n\frac{(-1)^{k-i}}{n!}
\binq{n}{i}
\sum_{j=0}^i\binm{i}{j}   
\lam^{i-j}(ky)^j\\
&=\sum_{m=0}^nx^m
\sum_{i=0}^n\frac{(-1)^{m-i}}{n!}
\binq{n}{i}
\sum_{j=0}^i\binm{i}{j}   
\lam^{i-j}y^j
\sum_{k=0}^m
(-1)^{m-k}\binm{m}{k}k^j.}
Writing the rightmost sum in terms of Stirling number of the second kind,
we get a triple sum expression
\pq{Q_n=\sum_{m=0}^nx^m
\sum_{i=0}^n(-1)^{m-i}\frac{m!}{n!}
\binq{n}{i}\sum_{j=0}^i\binm{i}{j}   
\bing{j}{m}\lam^{i-j}y^j.}
Taking into account that $\bing{j}{m}=0$ for $0\le j<m$
and $\binm{i}{j}=0$ for $0\le i<j$, we can restate
the last formula as follows.
\begin{prop}[Triple sum formula]\label{Q3s} 
\[Q_n=\sum_{m=0}^n\sum_{i=m}^n\sum_{j=m}^i(-1)^{m-i}\frac{m!}{n!}
\binq{n}{i}\binm{i}{j}\bing{j}{m}x^m\lam^{i-j}y^j.\]
\end{prop}

By comparing the coefficients of $x^m$ between the two expressions 
in Proposition \ref{Q2s} and Proposition \ref{Q3s}, we derive 
the following summation formula.
\begin{thm}[Convolution formula]\label{thm=S}
\[\sum_{i=m}^n\sum_{j=m}^i(-1)^{n+i}
\binq{n}{i}\binm{i}{j}\bing{j}{m}\lam^{i-j}y^j
=\frac{n!}{m!}\sum_{k=0}^m
(-1)^{m+k}\binm{m}{k}\binm{\lam+ky}{n}.\]
\end{thm}

\section{Applications to Convolution Formulae}
By specifying $y$ to particular values, we shall derive 
from Theorem~\ref{thm=S} two pairs of combinatorial 
identities, that extend, with an extra parameter 
``$\lam$", the related known convolution formulae 
on Stirling numbers.


The first pair of identities are given in the following theorem.
\begin{thm}[$1\le m\le n$]\label{thm=1}
\pnq{\emph{(a)}\quad&
\sum_{i=m}^n\sum_{j=m}^i(-1)^{n+i}
\binq{n}{i}\binm{i}{j}\bing{j}{m}\lam^{i-j}
=\frac{n!}{m!}\binm{\lam}{n-m},\\
\emph{(b)}\quad&
\sum_{i=m}^n\sum_{j=m}^i(-1)^{i+j}
\binq{n}{i}\binm{i}{j}\bing{j}{m}\lam^{i-j}
=\frac{n!}{m!}\binm{\lam-m}{n-m}(-1)^{m+n}.}
\end{thm}

They contain the two well--known results 
below as special cases of $\lam=0$.
\begin{corl}[$1\le m\le n$]
\pnq{\emph{(a)}\quad&
\sum_{k=m}^n(-1)^{n-k}
\binq{n}{k}\bing{k}{m}
=\begin{cases}
1,&m=n;\\
0,&m\not=n;
\end{cases}\\
\emph{(b)}\quad&
\sum_{k=m}^n\binq{n}{k}\bing{k}{m}
=\frac{n!}{m!}\binm{n-1}{m-1}.
\qquad\fbox{\emph{Lah number}}}
\end{corl}

\emph{Proof of Theorem~\ref{thm=1}.} \
Letting $y=1$ in Theorem~\ref{thm=S}, the corresponding
equality becomes 
\[\sum_{i=m}^n\sum_{j=m}^i(-1)^{n+i}
\binq{n}{i}\binm{i}{j}\bing{j}{m}\lam^{i-j}
=\frac{n!}{m!}\sum_{k=0}^m
(-1)^{m+k}\binm{m}{k}\binm{\lam+k}{n}.\]
Then the first formula ``(a)" follows by evaluating 
the above sum on the right
\pnq{\sum_{k=0}^m(-1)^{m+k}
\binm{m}{k}\binm{\lam+k}{n}
&=\sum_{k=0}^m(-1)^{m+k}
\binm{m}{k}\sum_{\ell=0}^n
\binm{\lam}{n-\ell}\binm{k}{\ell}\\
&=\sum_{\ell=0}^n\binm{\lam}{n-\ell}
\binm{m}{\ell}\sum_{k=\ell}^m(-1)^{m+k}
\binm{m-\ell}{k-\ell}\\
&=\sum_{\ell=0}^n\binm{\lam}{n-\ell}
\binm{m}{\ell}(1-1)^{m-\ell}
=\binm{\lam}{n-m},}
where the Chu--Vandemonde convolution formula 
and the binomial theorem have been utilized.

Alternatively for $y=-1$, the corresponding equality 
in Theorem~\ref{thm=S} reads as 
\[\sum_{i=m}^n\sum_{j=m}^i(-1)^{i+j}
\binq{n}{i}\binm{i}{j}\bing{j}{m}\lam^{i-j}
=(-1)^{m-n}\frac{n!}{m!}\sum_{k=0}^m
(-1)^{k}\binm{m}{k}\binm{\lam-k}{n}.\]
Then we can analogously evaluate the sum on the right as follows:
\pnq{\sum_{k=0}^m
(-1)^{k}\binm{m}{k}\binm{\lam-k}{n}
&=\sum_{k=0}^m
(-1)^{k}\binm{m}{k}
\sum_{\ell=0}^n
\binm{\lam-m}{n-\ell}\binm{m-k}{\ell}\\
&=\sum_{\ell=0}^n
\binm{\lam-m}{n-\ell}\binm{m}{\ell}
\sum_{k=0}^{m-\ell}
(-1)^{k}\binm{m-\ell}{k}\\
&=\sum_{\ell=0}^n
\binm{\lam-m}{n-\ell}\binm{m}{\ell}
(1-1)^{m-\ell}
=\binm{\lam-m}{n-m}.}
This proves the second formula ``(b)" in the theorem.\qed

We can further show another pair of identities. 
\begin{thm}[$1\le m\le n$]\label{thm=2}
\pnq{\emph{(a)}\quad
\sum_{i=m}^n&\sum_{j=m}^i(-1)^{n+i}
\binq{n}{i}\binm{i}{j}\bing{j}{m}\lam^{i-j}2^{n-j}\\
=~&\frac{n!}{m!}\sum_{\ell=m}^n\frac{m}{\ell}
\binm{\lam}{n-\ell}\binm{-\ell}{\ell-m}
2^{m+n-2\ell},\\
\emph{(b)}\quad
\sum_{i=m}^n&\sum_{j=m}^i(-1)^{m+i}
\binq{n}{i}\binm{i}{j}\bing{j}{m}\lam^{i-j}2^{j-m}\\
=~&\frac{n!}{m!}\sum_{\ell=m}^n
(-1)^{m+n}\binm{\lam}{n-\ell}
\binm{m}{\ell-m}2^{m-\ell}.}
\end{thm}

\emph{Proof of Theorem~\ref{thm=2}.} \
For $y=\frac12$, the equality in Theorem~\ref{thm=S} reduces to
\[\sum_{i=m}^n\sum_{j=m}^i(-1)^{n+i}
\binq{n}{i}\binm{i}{j}\bing{j}{m}\lam^{i-j}2^{n-j}
=2^n\frac{n!}{m!}\sum_{k=0}^m
(-1)^{m+k}\binm{m}{k}\binm{\lam+\frac{k}2}{n}.\]
In this case, the first formula ``(a)" is confirmed 
by first reformulating the above sum on the right
\pnq{\sum_{k=0}^m(-1)^{m+k}\binm{m}{k}\binm{\lam+\frac{k}2}{n}
&=\sum_{k=0}^m(-1)^{m+k}\binm{m}{k}
\sum_{\ell=0}^n\binm{\lam}{n-\ell}\binm{\frac{k}2}{\ell}\\
&=\sum_{\ell=0}^n\binm{\lam}{n-\ell}
\sum_{k=0}^m(-1)^{m+k}\binm{m}{k}\binm{\frac{k}2}{\ell}}
and then evaluating the inner sum by the binomial 
formula (cf.~\citu{gould}{Equation~3.164})
\[\sum_{k=0}^m
(-1)^{m-k}\binm{m}{k}\binm{\frac{k}2}{\ell}
=\frac{m}{\ell}\binm{-\ell}{\ell-m}2^{m-2\ell}.\]

When $y=2$, the corresponding equality in Theorem~\ref{thm=S} can be stated as 
\[\sum_{i=m}^n\sum_{j=m}^i(-1)^{m+i}
\binq{n}{i}\binm{i}{j}\bing{j}{m}\lam^{i-j}2^{j-m}
=2^{-m}\frac{n!}{m!}\sum_{k=0}^m
(-1)^{n+k}\binm{m}{k}\binm{\lam+2k}{n}.\]
The above sum on the right can be rewritten as
\pnq{\sum_{k=0}^m(-1)^{n+k}\binm{m}{k}\binm{\lam+2k}{n}
&=\sum_{k=0}^m(-1)^{n+k}\binm{m}{k}
\sum_{\ell=0}^n\binm{\lam}{n-\ell}\binm{2k}{\ell}\\
&=\sum_{\ell=0}^n(-1)^{n}\binm{\lam}{n-\ell}
\sum_{k=0}^m(-1)^{k}\binm{m}{k}\binm{2k}{\ell}.}

Denote by $[x^n]f(x)$ the coefficient of $x^n$ 
in the formal power series $f(x)$. Observe that
the last inner sum can be expressed and then evaluated as  
\pnq{\sum_{k=0}^m(-1)^{k}\binm{m}{k}\binm{2k}{\ell}
&=[x^{\ell}]\sum_{k=0}^m(-1)^{k}\binm{m}{k}(1+x)^{2k}\\
&=[x^{\ell}]\Big\{1-(1+x)^2\Big\}^m\\
&=(-1)^m[x^{\ell-m}](2+x)^m\\
&=(-1)^m\binm{m}{\ell-m}2^{2m-\ell}.}
After substitutions, the second formula ``(b)" 
in the theorem is done.\qed

Finally letting $\lam=0$ in Theorem~\ref{thm=2}, we recover the 
following two elegant identities due to Yang and Qiao~\cito{yang}, 
who discovered them by employing the Riordan array and expressed 
the results in terms of Bessel numbers~\cito{shapiro}.
\begin{corl}[$1\le m\le n$]
\pnq{\emph{(a)}\quad&
\sum_{k=m}^n\binq{n}{k}\bing{k}{m}(-2)^{n-k}
=\frac{(n-1)!}{(m-1)!}\binm{-n}{n-m}2^{m-n},\\
\emph{(b)}\quad&
\sum_{k=m}^n\binq{n}{k}\bing{k}{m}(-2)^{j-m}
=\frac{n!}{m!}\binm{m}{n-m}(-2)^{m-n}.}
\end{corl}


\thank{125mm}{The first author is partially supported, during this work, by 
the National Science foundation of China \emph{(Youth Grant No.11601543)}.}


\end{document}